\documentclass[10pt]{article}
\usepackage{latexsym, amsfonts, amsmath, amssymb, amscd, epsfig}

\numberwithin{equation}{section}
\newtheorem{thm}{Theorem}[section]
\newtheorem{defn}[thm]{Definition}

\newtheorem{lem}[thm]{Lemma}

\newtheorem{re}{Remark}[section]
\newenvironment{pf}{{\noindent \it \bf Proof:}}{{\hfill$\Box$}\\}

\begin{document}

\title{A regularity criterion in weak spaces to Boussinesq equations}
\author{Ravi P. Agarwal$\ ^{1}$, S. Gala$^{2,3}$ and Maria Alessandra Ragusa$%
^{3,4}$ \\
$^{1}${\small Department of Mathematics, Texas A \& M University-Kingsville,
kingsville, USA,}\\
$^{2}${\small University of Mostaganem, P. O. Box 270, Mostaganem 27000,
Algeria,}\\
$^{3}${\small Dipartimento di Matematica e Informatica, Universit\`{a} di
Catania,}\\
{\small Viale Andrea Doria, 6 95125 Catania - Italy,}\\
$^{4}${\small RUDN University, 6 Miklukho - Maklay St, Moscow, 117198,
Russia.}}
\date{}
\maketitle

\begin{abstract}
In this paper, we study regularity of weak solutions to the incompressible
Boussinesq equations in $\mathbb{R}^{3}\times (0,T)$. The main goal is to
establish the regularity criterion in terms of one velocity component and
the gradient of temperature in Lorentz spaces.
\end{abstract}

\section{Introduction}

In this paper we consider the following Cauchy problem for the
incompressible Boussinesq equations in $\mathbb{R}^{3}\times (0,T)$
\begin{equation}
\left\{
\begin{array}{c}
\partial _{t}u+\left( u\cdot \nabla \right) u-\Delta u+\nabla \pi =\theta
e_{3}, \\
\partial _{t}\theta -\Delta \theta +(u\cdot \nabla )\theta =0, \\
\nabla \cdot u=\nabla \cdot b=0, \\
u(x,0)=u_{0}(x),\text{ }\theta (x,0)=\theta _{0}(x),%
\end{array}%
\right.   \label{eq1}
\end{equation}%
where $u=(u_{1}(x,t),u_{2}(x,t),u_{3}(x,t))$ denotes the unknown velocity
vector, $\theta =(\theta _{1}(x,t),\theta _{2}(x,t),\theta _{3}(x,t))$ and $%
\pi =\pi (x,t)$ denote, respectively, the temperature and the hydrostatic
pressure. While $u_{0}$ and $\theta _{0}$ are the prescribed initial data
for the velocity and temperature with properties $\nabla \cdot u_{0}=0$.
Moreover, the term $\theta e_{3}$ represents buoyancy force on \ fluid
motion.

We would like to point out that the system (\ref{eq1}) at $\theta =0$
reduces to the incompressible Navier-Stokes equations, which has been
greatly analyzed. From the viewpoint of the model, therefore, Navier-Stokes
flow is viewed as the \ flow of a simplified Boussinesq equation.

Besides their physical applications, the Boussinesq equations are also
mathematically significant. Fundamental mathematical issues such as the
global regularity of their solutions have generated extensive research, and
many interesting results have been obtained (see, for example, \cite{CN,
CKN, FO, FZ, GGR, G1, GR, IM, QYW, QDY, X, XZZ, Ye3, Ye} and references
therein).

On the other hand, it is desirable to show the regularity of the weak
solutions if some partial components of the velocity satisfy certain growth
conditions. For the 3D Navier-Stokes equations, there are many results to
show such regularity of weak solutions in terms of partial components of the
velocity $u$ (see, for example, \cite{CT0, CT1, CZ, FQ1, FQ2, He, JZ1, PP,
Zhe} and the references cited therein). It is obvious that, for the
assumptions of all regularity criteria, it need that every components of the
velocity field must satisfies the same assumptions, and it don't make any
difference between the components of the velocity field. As pointed out by
Neustupa and Penel \cite{NP}, it is interesting to know how to effect the
regularity of the velocity field by the regularity of only one component of
the velocity field. In particular, Zhou \cite{Z1} showed that the solution
is regular if one component of the velocity, for example, $u_{3}$ satisfies
\begin{equation}
u_{3}\in L^{p}(0,T;L^{q}(\mathbb{R}^{3}))\text{ \ \ with \ }\frac{2}{p}+%
\frac{3}{q}\leq \frac{1}{2},\text{ \ \ }6<q\leq \infty .  \label{eq02}
\end{equation}%
Condition \ref{eq02}) can be replaced respectively by the one
\begin{equation}
u_{3}\in L^{p}(0,T;L^{q}(\mathbb{R}^{3}))\text{ \ \ with \ }\frac{2}{p}+%
\frac{3}{q}\leq \frac{5}{8},\text{ \ \ }\frac{24}{5}<q\leq \infty ,
\label{eq002}
\end{equation}%
(see Kukavica and Ziane \cite{KZ}). Later, Cao and Titi \cite{CT1} showed
the regularity of weak solution to the Navier-Stokes equations under the
assumption
\begin{equation}
u_{3}\in L^{p}(0,T;L^{q}(\mathbb{R}^{3}))\text{ \ \ with \ }\frac{2}{p}+%
\frac{3}{q}=\frac{2}{3}+\frac{2}{3q},\text{ \ \ }q>\frac{7}{2}.  \label{eq03}
\end{equation}%
Motivated by the above work, Zhou and Pokorn\'{y} \cite{ZP1} showed the
following regularity condition%
\begin{equation}
u_{3}\in L^{p}(0,T;L^{q}(\mathbb{R}^{3}))\text{ \ \ with \ }\frac{2}{p}+%
\frac{3}{q}=\frac{3}{4}+\frac{1}{2q},\text{ \ \ }q>\frac{10}{3},  \label{eq2}
\end{equation}%
while the limiting case $u_{3}\in L^{\infty }(0,T;L^{\frac{10}{3}}(\mathbb{R}%
^{3}))$ was covered in \cite{JZ1}. For many other result works, especially
the regularity criteria involving only one velocity component, or its
gradient, with no intention to be complete, one can consult \cite{Ye1, Ye2}
and references therein. However, the conditions (\ref{eq02})-(\ref{eq2}) are
quite strong comparing with the condition of Serrin's regularity criterion :%
\begin{equation}
u\in L^{p}(0,T;L^{q}(\mathbb{R}^{3}))\text{ \ \ with \ }\frac{2}{p}+\frac{3}{%
q}\leq 1,\text{ \ \ }3<q\leq \infty ,  \label{eq1.7}
\end{equation}%
and do not imply the invariance under the above scaling transformation.
Therefore, it is of interest in showing regularity by imposing Serrin's
condition (\ref{eq1.7}) with respect to the one component of the velocity
field.

Similar to the research of the 3D Navier-Stokes equations, authors are
interested in the regularity criterion of (\ref{eq1}) by reducing to some
the components of $u$. There are many other or similar results on the
hydrodynamical systems modeling the flow of nematic liquid crystal material,
Boussinesq equations and MHD equations (see e.g. \cite{BGR, Q} and the
reference therein).

Motivated by the reference mentioned above, the purpose of the present paper
is to give a further observation on the global regularity of the solution
for system (\ref{eq1}) and to extend the regularity of weak solutions to the
Boussinesq equations (\ref{eq1}) in terms of one velocity component and the
gradient of the temperature.

\section{Notations and main result.}

Before stating our result, we introduce some notations and function spaces.
These spaces can be found in many literatures and papers. For the functional
space, $L^{p}(\mathbb{R}^{3})$ denotes the usual Lebesgue space of
real-valued functions with norm $\left\Vert \cdot \right\Vert _{L^{p}}:$%
\begin{equation*}
\left\Vert f\right\Vert _{L^{p}}=\left\{
\begin{array}{c}
\left( \int_{\mathbb{R}^{3}}\left\vert f(x)\right\vert ^{p}dx\right) ^{\frac{%
1}{p}},\text{ \ \ for \ \ }1\leq p<\infty , \\
\underset{x\in \mathbb{R}^{3}}{ess\sup }\left\vert f(x)\right\vert ,\text{ \
\ for \ \ }p=\infty .%
\end{array}%
\right.
\end{equation*}%
On the other hand, the usual Sobolev space of order $m$ is defined by
\begin{equation*}
H^{m}(\mathbb{R}^{3})=\left\{ u\in L^{2}(\mathbb{R}^{3}):\nabla ^{m}u\in
L^{2}(\mathbb{R}^{3})\right\}
\end{equation*}%
with the norm%
\begin{equation*}
\left\Vert u\right\Vert _{H^{m}}=\left( \left\Vert u\right\Vert
_{L^{2}}^{2}+\left\Vert \nabla ^{m}u\right\Vert _{L^{2}}^{2}\right) ^{\frac{1%
}{2}}.
\end{equation*}%
To prove Theorem \ref{th1}, we use the theory of Lorentz spaces and
introduce the following notations. We define the non-increasing
rearrangement of $f$,%
\begin{equation*}
f^{\ast }(\lambda )=\inf \left\{ t>0:m_{f}(t)\leq \lambda \right\} ,\text{ \
\ for \ }\lambda >0,
\end{equation*}%
where $f$ is a measurable function on $\mathbb{R}^{3}$ and $m_{f}(t)$ is the
distribution function of $f$ which is defined by the Lebesgue measure of the
set $\left\{ x\in \mathbb{R}^{3}:\left\vert f(x)\right\vert >t\right\} $.
The Lorentz space $L^{p,q}((\mathbb{R}^{3})$ is defined by
\begin{equation*}
L^{p,q}=\left\{ f:\mathbb{R}^{3}\rightarrow \mathbb{R}\text{ measurable such
that }\left\Vert f\right\Vert _{L^{p,q}}<\infty \right\} \text{ \ with \ }%
1\leq p<\infty
\end{equation*}%
equipped with the quasi-norm
\begin{equation*}
\left\Vert f\right\Vert _{L^{p,q}}=\left( \frac{q}{p}\int_{0}^{\infty }(t^{%
\frac{1}{p}}f^{\ast }(t))^{q}\frac{dt}{t}\right) ^{\frac{1}{q}}\text{, \ \
if \ }1<q<\infty .
\end{equation*}%
Moreover, we define $f^{\ast \ast }$ by
\begin{equation*}
f^{\ast \ast }(\lambda )=\frac{1}{\lambda }\int_{0}^{\lambda }f^{\ast
}(\lambda ^{\prime })d\lambda ^{\prime },
\end{equation*}%
and Lorentz spaces $L^{p,\infty }(\mathbb{R}^{3})$ by%
\begin{equation*}
L^{p,\infty }(\mathbb{R}^{3})=\left\{ f\in \mathcal{S}^{\prime }(\mathbb{R}%
^{3}):\left\Vert f\right\Vert _{L^{p,\infty }}<\infty \right\} ,
\end{equation*}%
where
\begin{equation*}
\left\Vert f\right\Vert _{L^{p,\infty }}=\underset{\lambda \geq 0}{\sup }%
(\lambda ^{\frac{1}{p}}f^{\ast \ast }(\lambda )),
\end{equation*}%
for $1\leq p\leq \infty $. For details, we refer to \cite{BL} and \cite{Tri}.

From the definition of the Lorentz space, we can obtain the following
continuous embeddings :%
\begin{equation*}
L^{p}(\mathbb{R}^{3})=L^{p,p}(\mathbb{R}^{3})\hookrightarrow L^{p,q}(\mathbb{%
R}^{3})\hookrightarrow L^{p,\infty }(\mathbb{R}^{3}),\text{ \ \ }1\leq p\leq
q<\infty .
\end{equation*}%
In order to prove Theorem \ref{th1}, we recall the H\"{o}lder inequality in
the Lorentz spaces (see, e.g., O'Neil \cite{O}).

\begin{lem}
\label{lem1}Let $f\in L^{p_{2},q_{2}}(\mathbb{R}^{3})$ and $g\in
L^{p_{3},q_{3}}(\mathbb{R}^{3})$ with $1\leq p_{2},p_{3}\leq \infty $, $%
1\leq q_{2},q_{3}\leq \infty $.\ Then $fg\in L^{p_{1},q_{1}}(\mathbb{R}^{3})$
with
\begin{equation*}
\frac{1}{p_{1}}=\frac{1}{p_{2}}+\frac{1}{p_{3}},\frac{1}{q_{1}}=\frac{1}{%
q_{2}}+\frac{1}{q_{3}}
\end{equation*}%
and the H\"{o}lder inequality of Lorentz spaces
\begin{equation*}
\left\Vert fg\right\Vert _{L^{p_{1},q_{1}}}\leq C\left\Vert f\right\Vert
_{L^{p_{2},q_{2}}}\left\Vert g\right\Vert _{L^{p_{3},q_{3}}},
\end{equation*}%
holds true for a positive constant $C$.
\end{lem}

The following result plays an important role in the proof of our theorem,
the so-called Gagliardo-Nirenberg inequality in Lorentz spaces, its proof
can be founded in \cite{H}.

\begin{lem}
\label{lem2}Let $f\in L^{p,q}(\mathbb{R}^{3})$ with $1\leq
p,q,p_{4},q_{4},p_{5},q_{5}\leq \infty $. Then the Gagliardo-Nirenberg
inequality of Lorentz spaces%
\begin{equation*}
\left\Vert f\right\Vert _{L^{p,q}}\leq C\left\Vert f\right\Vert
_{L^{p_{4},q_{4}}}^{\theta }\left\Vert f\right\Vert
_{L^{p_{5},q_{5}}}^{1-\theta }
\end{equation*}%
holds for a positive constant $C$ and%
\begin{equation*}
\frac{1}{p}=\frac{\theta }{p_{4}}+\frac{1-\theta }{p_{5}},\text{ \ \ }\frac{1%
}{q}=\frac{\theta }{q_{4}}+\frac{1-\theta }{q_{5}},\text{ \ }\theta \in
(0,1).
\end{equation*}
\end{lem}

Now we give the definition of weak solution.

\begin{defn}
Let $T>0$, $(u_{0},\theta _{0})\in L^{2}(\mathbb{R}^{3})$ with $\nabla \cdot
u_{0}=0$ in the sense of distributions. A measurable function $%
(u(x,t),\theta (x,t))$ is called a weak solution to the Boussinesq equations
(\ref{eq1}) on $[0,T]$ if the following conditions hold:

\begin{enumerate}
\item $(u(x,t),\theta (x,t))\in L^{\infty }(0,T;L^{2}(\mathbb{R}^{3}))\cap
L^{2}(0,T;H^{1}(\mathbb{R}^{3}));$

\item system (\ref{eq1}) is satisfied in the sense of distributions;

\item the energy inequality, that is,
\begin{equation*}
\left\Vert u(\cdot ,t)\right\Vert _{L^{2}}^{2}+\left\Vert \theta (\cdot
,t)\right\Vert _{L^{2}}^{2}+2\int_{0}^{t}\left\Vert \nabla u(\tau
)\right\Vert _{L^{2}}^{2}d\tau +2\int_{0}^{t}\left\Vert \nabla \theta (\tau
)\right\Vert _{L^{2}}^{2}d\tau \leq \left\Vert u_{0}\right\Vert
_{L^{2}}^{2}+\left\Vert b_{0}\right\Vert _{L^{2}}^{2}+\left\Vert \theta
_{0}\right\Vert _{L^{2}}^{2}.
\end{equation*}
\end{enumerate}
\end{defn}

By a strong solution, we mean that a weak solution u of the Navier-Stokes
equations (\ref{eq1}) satisfies%
\begin{equation*}
(u(x,t),\theta (x,t))\in L^{\infty }(0,T;H^{1}(\mathbb{R}^{3}))\cap
L^{2}(0,T;H^{2}(\mathbb{R}^{3})).
\end{equation*}%
It is well known that the strong solution is regular and unique.

Our main result is stated as following :

\begin{thm}
\label{th1}Let $(u_{0},\theta _{0})\in L^{2}(\mathbb{R}^{3})$ with $\nabla
\cdot u_{0}=0$ in the sense of distributions. Assume that $(u,\theta )$ is a
weak solution to system (\ref{eq1}).\ If $u_{3}$ and $\nabla \theta $
satisfy the following conditions%
\begin{equation}
\left\{
\begin{array}{c}
u_{3}\in L^{\frac{30\alpha }{7\alpha -45}}(0,T;L^{\alpha ,\infty }(\mathbb{R}%
^{3})),\text{ \ with \ }\frac{45}{7}\leq \alpha \leq \infty , \\
\nabla \theta \in L^{\frac{2\beta }{2\beta -3}}(0,T;L^{\beta ,\infty }(%
\mathbb{R}^{3})),\text{ \ \ with \ \ }\frac{3}{2}<\beta \leq \infty ,%
\end{array}%
\right.   \label{eq16}
\end{equation}%
then the solution $\left( u,\theta \right) $ is regular on $(0,T]$.
\end{thm}

\begin{re}
If $\theta =0$, it is clear that theorem \ref{th1} improves the earlier
results of \cite{JZ1, ZP1} for 3D Navier-Stokes equations and extend the
regularity criterion (\ref{eq2}) from Lebesgue space $L^{\alpha }$ to
Lorentz space $L^{\alpha ,\infty }$.
\end{re}

\begin{re}
This result proves a new regularity criterion for weak solutions to the
Cauchy problem of the 3D Boussinesq equations via one component of the
velocity field and the gradient of the tempearture in the framework of the
Lorentz spaces. This result reveals that the one component of the velocity
field plays a dominant role in regularity theory of the Boussinesq equations.
\end{re}

\section{Proof of the main result.}

In this section, under the assumptions of the Theorem \ref{th1}, we prove
our main result. Before proving our result, we recall the following
muliplicative Sobolev imbedding inequality in the whole space $\mathbb{R}%
^{3} $ (see, for example \cite{CT1}) :
\begin{equation}
\left\Vert f\right\Vert _{L^{6}}\leq C\left\Vert \nabla _{h}f\right\Vert
_{L^{2}}^{\frac{2}{3}}\left\Vert \partial _{3}f\right\Vert _{L^{2}}^{\frac{1%
}{3}},  \label{eq8}
\end{equation}%
where $\nabla _{h}=(\partial _{x_{1}},\partial _{x_{2}})$ is the horizontal
gradient operator. We are now give the proof of our main theorem.

\begin{pf}
To prove our result, it suffices to show that for any fixed $T>T^{\ast }$,
there holds%
\begin{equation*}
\underset{0\leq t\leq T^{\ast }}{\sup }(\left\Vert \nabla u(t)\right\Vert
_{L^{2}}^{2}+\left\Vert \nabla \theta (t)\right\Vert _{L^{2}}^{2})\leq C_{T},
\end{equation*}%
where $T^{\ast }$, which denotes the maximal existence time of a strong
solution and $C_{T}$ is an absolute constant which only depends on $T,u_{0}$
and $\theta _{0}$.

The method of our proof is based on two major parts. The first one
establishes the bounds of $(\left\Vert \nabla _{h}u\right\Vert
_{L^{2}}^{2}+\left\Vert \nabla _{h}\theta (t)\right\Vert _{L^{2}}^{2})$,
while the second gives the bounds of the $H^{1}-$norm of velocity $u$ and
temperature $\theta $\ in terms of the results of part one.

Taking the inner product of (\ref{eq1})$_{1}$ with $-\Delta _{h}u$, (\ref%
{eq1})$_{2}$ with $-\Delta _{h}\theta $ in $L^{2}(\mathbb{R}^{3})$,
respectively, then adding the three resulting equations together, we obtain
after integrating by parts that
\begin{eqnarray}
&&\frac{1}{2}\frac{d}{dt}(\left\Vert \nabla _{h}u\right\Vert
_{L^{2}}^{2}+\left\Vert \nabla _{h}\theta \right\Vert
_{L^{2}}^{2})+\left\Vert \nabla \nabla _{h}u\right\Vert
_{L^{2}}^{2}+\left\Vert \nabla \nabla _{h}\theta \right\Vert _{L^{2}}^{2}
\notag \\
&=&\int_{\mathbb{R}^{3}}(u\cdot \nabla )u\cdot \Delta _{h}udx+\int_{\mathbb{R%
}^{3}}(u\cdot \nabla )\theta \cdot \Delta _{h}\theta -\int_{\mathbb{R}%
^{3}}(\theta e_{3})\cdot \Delta _{h}udx  \notag \\
&=&I_{1}+I_{2}+I_{3}.  \label{eq5}
\end{eqnarray}%
where $\Delta _{h}=\partial _{x_{1}}^{2}+\partial _{x_{2}}^{2}$ is the
horizontal Laplacian. For the notational simplicity, we set%
\begin{eqnarray*}
\mathcal{L}^{2}(t) &=&\underset{\tau \in \lbrack \Gamma ,t]}{\sup }%
(\left\Vert \nabla _{h}u(\tau )\right\Vert _{L^{2}}^{2}+\left\Vert \nabla
_{h}\theta \right\Vert _{L^{2}}^{2})+\int_{\Gamma }^{t}(\left\Vert \nabla
\nabla _{h}u(\tau )\right\Vert _{L^{2}}^{2}+\left\Vert \nabla \nabla
_{h}\theta \right\Vert _{L^{2}}^{2})d\tau , \\
\mathcal{J}^{2}(t) &=&\underset{\tau \in \lbrack \Gamma ,t]}{\sup }%
(\left\Vert \nabla u(\tau )\right\Vert _{L^{2}}^{2}+\left\Vert \nabla \theta
(\tau )\right\Vert _{L^{2}}^{2})+\int_{\Gamma }^{t}(\left\Vert \Delta u(\tau
)\right\Vert _{L^{2}}^{2}+\left\Vert \Delta \theta \right\Vert
_{L^{2}}^{2})d\tau ,
\end{eqnarray*}%
for $t\in \lbrack \Gamma ,T^{\ast })$. In view of (\ref{eq16}), we choose $%
\epsilon ,\eta >0$ to be precisely determined subsequently and then select $%
\Gamma <T^{\ast }$ sufficiently close to $T^{\ast }$ such that for all $%
\Gamma \leq t<T^{\ast }$,%
\begin{equation}
\int_{\Gamma }^{t}(\left\Vert \nabla u(\tau )\right\Vert
_{L^{2}}^{2}+\left\Vert \nabla \theta (\tau )\right\Vert _{L^{2}}^{2})d\tau
\leq \epsilon \ll 1\text{,}  \label{eq19}
\end{equation}%
and%
\begin{equation}
\int_{\Gamma }^{t}\left\Vert \nabla \theta (\tau )\right\Vert _{L^{\beta }}^{%
\frac{2\beta }{2\beta -3}}d\tau \leq \eta \ll 1.  \label{eq23}
\end{equation}%
Integrating by parts and using the divergence-free condition, it is clear
that (see e.g. \cite{Zhe})
\begin{equation}
I_{1}=\int_{\mathbb{R}^{3}}(u\cdot \nabla )u\cdot \Delta _{h}udx\leq \int_{%
\mathbb{R}^{3}}\left\vert u_{3}\right\vert \left\vert \nabla u\right\vert
\left\vert \nabla \nabla _{h}u\right\vert dx  \label{eq3}
\end{equation}%
By appealing to Lemma \ref{lem1}, (\ref{eq8}), and the Young inequality, it
follows that
\begin{eqnarray*}
I_{1} &\leq &C\left\Vert u_{3}\right\Vert _{L^{\alpha ,\infty }}\left\Vert
\nabla u\right\Vert _{L^{\frac{2\alpha }{\alpha -2},2}}\left\Vert \nabla
\nabla _{h}u\right\Vert _{L^{2}} \\
&\leq &C\left\Vert u_{3}\right\Vert _{L^{\alpha ,\infty }}\left\Vert \nabla
u\right\Vert _{L^{2}}^{1-\frac{3}{\alpha }}\left\Vert \nabla u\right\Vert
_{L^{6}}^{\frac{3}{\alpha }}\left\Vert \nabla \nabla _{h}u\right\Vert
_{L^{2}} \\
&\leq &C\left\Vert u_{3}\right\Vert _{L^{\alpha ,\infty }}\left\Vert \nabla
u\right\Vert _{L^{2}}^{1-\frac{3}{\alpha }}\left\Vert \Delta u\right\Vert
_{L^{2}}^{\frac{1}{\alpha }}\left\Vert \nabla \nabla _{h}u\right\Vert
_{L^{2}}^{1+\frac{2}{\alpha }} \\
&\leq &C\left\Vert u_{3}\right\Vert _{L^{\alpha ,\infty }}^{\frac{2\alpha }{%
\alpha -2}}\left\Vert \nabla u\right\Vert _{L^{2}}^{2-\frac{2}{\alpha -2}%
}\left\Vert \Delta u\right\Vert _{L^{2}}^{\frac{2}{\alpha -2}}+\frac{1}{4}%
\left\Vert \nabla \nabla _{h}u\right\Vert _{L^{2}}^{2},
\end{eqnarray*}%
where we have used the following Gagliardo-Nirenberg inequality in Lorentz
spaces :
\begin{equation*}
\left\Vert \phi \right\Vert _{L^{\frac{2s}{s-2},2}}\leq C\left\Vert \phi
\right\Vert _{L^{2}}^{1-\frac{3}{s}}\left\Vert \nabla \phi \right\Vert
_{L^{6}}^{\frac{3}{s}}\text{ \ \ with \ \ }3<s\leq \infty .
\end{equation*}%
To estimate the term $I_{2}$ of (\ref{eq5}), first observe that by applying
integration by parts and $\nabla \cdot u=0$, we derive
\begin{eqnarray*}
I_{2} &=&\int_{\mathbb{R}^{3}}(u\cdot \nabla )\theta \cdot \Delta _{h}\theta
dx=\sum\limits_{i,j=1}^{3}\sum\limits_{k=1}^{2}\int_{\mathbb{R}%
^{3}}u_{i}\partial _{i}\theta _{j}\partial _{kk}^{2}\theta _{j}dx \\
&=&-\sum\limits_{i,j=1}^{3}\sum\limits_{k=1}^{2}\int_{\mathbb{R}%
^{3}}\partial _{k}u_{i}\partial _{i}\theta _{j}\partial _{k}\theta
_{j}dx-\sum\limits_{i,j=1}^{3}\sum\limits_{k=1}^{2}\int_{\mathbb{R}%
^{3}}u_{i}\partial _{k}\partial _{i}\theta _{j}\partial _{k}\theta _{j}dx \\
&=&-\sum\limits_{i,j=1}^{3}\sum\limits_{k=1}^{2}\int_{\mathbb{R}%
^{3}}\partial _{k}u_{i}\partial _{i}\theta _{j}\partial _{k}\theta _{j}dx,
\end{eqnarray*}%
where we have used
\begin{eqnarray*}
-\sum\limits_{i,j=1}^{3}\sum\limits_{k=1}^{2}\int_{\mathbb{R}%
^{3}}u_{i}\partial _{k}\partial _{i}\theta _{j}\partial _{k}\theta _{j}dx
&=&\sum\limits_{i,j=1}^{3}\sum\limits_{k=1}^{2}\int_{\mathbb{R}^{3}}\partial
_{i}u_{i}(\partial _{k}\theta
_{j})^{2}dx+\sum\limits_{i,j=1}^{3}\sum\limits_{k=1}^{2}\int_{\mathbb{R}%
^{3}}u_{i}\partial _{i}\partial _{k}\theta _{j}\partial _{k}\theta _{j}dx \\
&=&\sum\limits_{i,j=1}^{3}\sum\limits_{k=1}^{2}\int_{\mathbb{R}%
^{3}}u_{i}\partial _{i}\partial _{k}\theta _{j}\partial _{k}\theta _{j}dx,
\end{eqnarray*}%
so that
\begin{equation*}
\sum\limits_{i,j=1}^{3}\sum\limits_{k=1}^{2}\int_{\mathbb{R}%
^{3}}u_{i}\partial _{k}\partial _{i}\theta _{j}\partial _{k}\theta _{j}dx=0.
\end{equation*}%
It follows from H\"{o}lder's inequality, (\ref{eq8}) and Young's inequality
that%
\begin{eqnarray*}
I_{2} &=&-\sum\limits_{i,j=1}^{3}\sum\limits_{k=1}^{2}\int_{\mathbb{R}%
^{3}}\partial _{k}u_{i}\partial _{i}\theta _{j}\partial _{k}\theta
_{j}dx\leq \int_{\mathbb{R}^{3}}\left\vert \nabla _{h}u\right\vert
\left\vert \nabla \theta \right\vert \left\vert \nabla _{h}\theta
\right\vert dx \\
&\leq &\frac{1}{2}\int_{\mathbb{R}^{3}}\left( \left\vert \nabla
_{h}u\right\vert ^{2}+\left\vert \nabla _{h}\theta \right\vert ^{2}\right)
\left\vert \nabla \theta \right\vert dx \\
&\leq &C\left\Vert \nabla \theta \right\Vert _{L^{\beta ,\infty
}}(\left\Vert \nabla _{h}u\right\Vert _{L^{\frac{2\beta }{\beta -2}%
,2}}\left\Vert \nabla _{h}u\right\Vert _{L^{2}}+\left\Vert \nabla _{h}\theta
\right\Vert _{L^{\frac{2\beta }{\beta -2},2}}\left\Vert \nabla _{h}\theta
\right\Vert _{L^{2}}) \\
&\leq &C\left\Vert \nabla \theta \right\Vert _{L^{\beta ,\infty
}}(\left\Vert \nabla _{h}u\right\Vert _{L^{2}}^{2-\frac{3}{\beta }%
}\left\Vert \nabla \nabla _{h}u\right\Vert _{L^{2}}^{\frac{3}{s}}+\left\Vert
\nabla _{h}\theta \right\Vert _{L^{2}}^{2-\frac{3}{\beta }}\left\Vert \nabla
\nabla _{h}\theta \right\Vert _{L^{2}}^{\frac{3}{s}}) \\
&\leq &C\left\Vert \nabla \theta \right\Vert _{L^{\beta ,\infty }}^{\frac{%
2\beta }{2\beta -3}}(\left\Vert \nabla _{h}u\right\Vert
_{L^{2}}^{2}+\left\Vert \nabla _{h}\theta \right\Vert _{L^{2}}^{2})+\frac{1}{%
2}\left\Vert \nabla \nabla _{h}\theta \right\Vert _{L^{2}}^{2}+\frac{1}{4}%
\left\Vert \nabla \nabla _{h}u\right\Vert _{L^{2}}^{2}.
\end{eqnarray*}%
Finally, we we want to estimate $I_{3}$. It follows from integration by
parts and Cauchy inequality that%
\begin{eqnarray*}
I_{3} &=&-\int_{\mathbb{R}^{3}}(\theta e_{3})\cdot \Delta _{h}udx=\int_{%
\mathbb{R}^{3}}\nabla _{h}(\theta e_{3})\cdot \nabla _{h}udx \\
&\leq &2\left\Vert \nabla _{h}u\right\Vert _{L^{2}}^{2}+2\left\Vert \nabla
_{h}\theta \right\Vert _{L^{2}}^{2}.
\end{eqnarray*}%
Inserting all the estimates into (\ref{eq3}), Gronwall's type argument using
\begin{equation*}
1\leq \underset{\lambda \in \lbrack \Gamma ,\tau ]}{\sup }\exp \left(
c\int_{\lambda }^{\tau }\left\Vert \nabla \theta (\varphi )\right\Vert
_{L^{\beta ,\infty }}^{\frac{2\beta }{2\beta -3}}d\varphi \right) \lesssim
\exp \left( c\int_{0}^{T^{\ast }}\left\Vert \nabla \theta (\varphi
)\right\Vert _{L^{\beta ,\infty }}^{\frac{2\beta }{2\beta -3}}d\varphi
\right) \lesssim 1,
\end{equation*}%
due to (\ref{eq16}) leads to, for every $\tau \in \lbrack \Gamma ,t]$%
\begin{eqnarray}
\mathcal{L}^{2}(t) &\leq &C+C\int_{\Gamma }^{t}\left\Vert u_{3}\right\Vert
_{L^{\alpha ,\infty }}^{\frac{2\alpha }{\alpha -2}}\left\Vert \nabla
u\right\Vert _{L^{2}}^{2-\frac{2}{\alpha -2}}\left\Vert \Delta u\right\Vert
_{L^{2}}^{\frac{2}{\alpha -2}}d\tau +C\int_{\Gamma }^{t}\left\Vert \nabla
\theta \right\Vert _{L^{\beta ,\infty }}^{\frac{2\beta }{2\beta -3}%
}(\left\Vert \nabla _{h}u\right\Vert _{L^{2}}^{2}+\left\Vert \nabla
_{h}\theta \right\Vert _{L^{2}}^{2})d\tau  \notag \\
&=&C+\mathcal{I}_{1}(t)+\mathcal{I}_{2}(t).  \label{eq4}
\end{eqnarray}%
Next, we analyze the right-hand side of (\ref{eq4}) one by one. First, due
to (\ref{eq19}) and the definition of $\mathcal{J}^{2}$, we have
\begin{eqnarray*}
\mathcal{I}_{1}(t) &\leq &C\left( \underset{\tau \in \lbrack \Gamma ,t]}{%
\sup }\left\Vert \nabla u(\tau )\right\Vert _{L^{2}}^{\frac{3}{2}-\frac{2}{%
\alpha -2}}\right) \int_{\Gamma }^{t}\left\Vert u_{3}(\tau )\right\Vert
_{L^{\alpha ,\infty }}^{\frac{2\alpha }{\alpha -2}}\left\Vert \nabla u(\tau
)\right\Vert _{L^{2}}^{\frac{1}{2}}\left\Vert \Delta u(\tau )\right\Vert
_{L^{2}}^{\frac{2}{\alpha -2}}d\tau \\
&\leq &C\mathcal{J}^{\frac{3}{2}-\frac{2}{\alpha -2}}(t)\left( \int_{\Gamma
}^{t}\left\Vert u_{3}(\tau )\right\Vert _{L^{\alpha ,\infty }}^{\frac{%
8\alpha }{3\alpha -10}}d\tau \right) ^{\frac{3}{4}-\frac{1}{\alpha -2}%
}\left( \int_{\Gamma }^{t}\left\Vert \nabla u(\tau )\right\Vert
_{L^{2}}^{2}d\tau \right) ^{\frac{1}{4}}\left( \int_{\Gamma }^{t}\left\Vert
\Delta u(\tau )\right\Vert _{L^{2}}^{2}d\tau \right) ^{\frac{1}{\alpha -2}}
\\
&\leq &C\mathcal{J}^{\frac{3}{2}-\frac{2}{\alpha -2}}(t)\left( \int_{\Gamma
}^{t}\left\Vert u_{3}(\tau )\right\Vert _{L^{\alpha ,\infty }}^{\frac{%
8\alpha }{3\alpha -10}}d\tau \right) ^{\frac{3}{4}-\frac{1}{\alpha -2}%
}\epsilon ^{\frac{1}{4}}\mathcal{J}^{\frac{2}{\alpha -2}}(t) \\
&=&C\epsilon ^{\frac{1}{4}}\mathcal{J}^{\frac{3}{2}}(t)\left( \int_{\Gamma
}^{t}\left\Vert u_{3}(\tau )\right\Vert _{L^{\alpha ,\infty }}^{\frac{%
8\alpha }{3\alpha -10}}d\tau \right) ^{\frac{3}{4}-\frac{1}{\alpha -2}}.
\end{eqnarray*}%
Finally, we deal with the term $\mathcal{I}_{2}(t)$. Applying H\"{o}lder and
Young inequalities, one has
\begin{eqnarray*}
\mathcal{I}_{2}(t) &\leq &C\underset{\tau \in \lbrack \Gamma ,t]}{\sup }%
(\left\Vert \nabla _{h}u(\tau )\right\Vert _{L^{2}}^{2}+\left\Vert \nabla
_{h}\theta (\tau )\right\Vert _{L^{2}}^{2})\int_{\Gamma }^{t}\left\Vert
\nabla \theta (\tau )\right\Vert _{L^{\beta ,\infty }}^{\frac{2\beta }{%
2\beta -3}}d\tau \\
&\leq &C\eta \mathcal{L}^{2}(t).
\end{eqnarray*}%
Hence, choosing $\eta $ small enough such that $C\eta <1$ and inserting the
above estimates of $\mathcal{I}_{1}(t)$ and $\mathcal{I}_{2}(t)$\ into (\ref%
{eq4}), we derive that for all $\Gamma \leq t<T^{\ast }:$%
\begin{eqnarray*}
\mathcal{L}^{2}(t) &\leq &C+C\epsilon ^{\frac{1}{4}}\mathcal{J}^{\frac{3}{2}%
}(t)\left( \int_{\Gamma }^{t}\left\Vert u_{3}(\tau )\right\Vert _{L^{\alpha
,\infty }}^{\frac{8\alpha }{3\alpha -10}}d\tau \right) ^{\frac{3\alpha -10}{%
4(\alpha -2)}} \\
&\leq &C+C\epsilon ^{\frac{1}{4}}\mathcal{J}^{\frac{3}{2}}(t)\left(
\int_{\Gamma }^{t}1+\left\Vert u_{3}(\tau )\right\Vert _{L^{\alpha }}^{\frac{%
30\alpha }{7\alpha -45}}d\tau \right) ^{\frac{3\alpha -10}{4(\alpha -2)}},
\end{eqnarray*}

which leads to%
\begin{equation}
\mathcal{L}^{2}(t)\leq C+C\epsilon ^{\frac{1}{4}}\mathcal{J}^{\frac{3}{2}%
}(t).  \label{eq25}
\end{equation}

Now, we will establish the bounds of $H^{1}-$norm of the velocity magnetic
field and micro-rotationel velocity. In order to do it, taking the inner
product of (\ref{eq1})$_{1}$ with $-\Delta u$, (\ref{eq1})$_{2}$ with $%
-\Delta b$ and (\ref{eq1})$_{3}$ with $-\Delta \theta $ in $L^{2}(\mathbb{R}%
^{3})$, respectively. Then, integration by parts gives the following
identity:
\begin{eqnarray*}
&&\frac{1}{2}\frac{d}{dt}(\left\Vert \nabla u\right\Vert
_{L^{2}}^{2}+\left\Vert \nabla \theta \right\Vert _{L^{2}}^{2})+\left\Vert
\Delta u\right\Vert _{L^{2}}^{2}+\left\Vert \Delta \theta \right\Vert
_{L^{2}}^{2} \\
&=&\int_{\mathbb{R}^{3}}(u\cdot \nabla )u\cdot \Delta udx+\int_{\mathbb{R}%
^{3}}(u\cdot \nabla )\theta \cdot \Delta \theta dx-\int_{\mathbb{R}%
^{3}}(\theta e_{3})\cdot \Delta udx.
\end{eqnarray*}%
Integrating by parts and using the divergence-free condition, one can easily
deduce that (see e.g. \cite{ZP1})
\begin{equation*}
\int_{\mathbb{R}^{3}}(u\cdot \nabla )u\cdot \Delta udx\leq C\int_{\mathbb{R}%
^{3}}\left\vert \nabla _{h}u\right\vert \left\vert \nabla u\right\vert
^{2}dx.
\end{equation*}%
We treat now the $\int_{\mathbb{R}^{3}}(u\cdot \nabla )\theta \cdot \Delta
\theta dx$-term. By integration by parts, we have
\begin{eqnarray}
\int_{\mathbb{R}^{3}}(u\cdot \nabla )\theta \cdot \Delta \theta dx
&=&-\sum\limits_{i=1}^{2}\sum\limits_{j,k=1}^{3}\int_{\mathbb{R}%
^{3}}\partial _{i}u_{k}\cdot \partial _{k}\theta _{j}\cdot \partial
_{i}\theta _{j}dx-\sum\limits_{j,k=1}^{3}\int_{\mathbb{R}^{3}}\partial
_{3}u_{k}\cdot \partial _{k}\theta _{j}\cdot \partial _{3}\theta _{j}dx
\notag \\
&=&-\sum\limits_{i=1}^{2}\sum\limits_{j,k=1}^{3}\int_{\mathbb{R}%
^{3}}\partial _{i}u_{k}\cdot \partial _{k}\theta _{j}\cdot \partial
_{i}\theta _{j}dx-\sum\limits_{k=1}^{2}\sum\limits_{j=1}^{3}\int_{\mathbb{R}%
^{3}}\partial _{3}u_{k}\cdot \partial _{k}\theta _{j}\cdot \partial
_{3}\theta _{j}dx  \notag \\
&&-\sum\limits_{j=1}^{3}\int_{\mathbb{R}^{3}}\partial _{3}u_{3}\cdot
\partial _{3}\theta _{j}\cdot \partial _{3}\theta _{j}dx  \notag \\
&=&-\sum\limits_{i=1}^{2}\sum\limits_{j,k=1}^{3}\int_{\mathbb{R}%
^{3}}\partial _{i}u_{k}\cdot \partial _{k}\theta _{j}\cdot \partial
_{i}\theta _{j}dx-\sum\limits_{k=1}^{2}\sum\limits_{j=1}^{3}\int_{\mathbb{R}%
^{3}}\partial _{3}u_{k}\cdot \partial _{k}\theta _{j}\cdot \partial
_{3}\theta _{j}dx  \notag \\
&&+\sum\limits_{i=1}^{2}\sum\limits_{j=1}^{3}\int_{\mathbb{R}^{3}}\partial
_{i}u_{i}\cdot \partial _{3}\theta _{j}\cdot \partial _{3}\theta _{j}dx
\notag \\
&=&R_{1}+R_{2}+R_{3}.  \label{eq3.3}
\end{eqnarray}%
Therefore, we have%
\begin{equation*}
\left\vert R_{1}+R_{3}\right\vert \leq \int_{\mathbb{R}^{3}}\left\vert
\nabla _{h}u\right\vert \left\vert \nabla \theta \right\vert ^{2}dx,
\end{equation*}%
and
\begin{equation*}
\left\vert R_{2}\right\vert \leq \int_{\mathbb{R}^{3}}\left\vert \nabla
u\right\vert \left\vert \nabla \theta \right\vert \left\vert \nabla
_{h}\theta \right\vert dx\leq \frac{1}{2}\int_{\mathbb{R}^{3}}\left\vert
\nabla _{h}\theta \right\vert (\left\vert \nabla u\right\vert
^{2}+\left\vert \nabla \theta \right\vert ^{2})dx,
\end{equation*}%
where the last inequality is obtained by using Cauchy inequality.

Putting all the inequalities above into (\ref{eq3.3}) yields%
\begin{equation*}
\int_{\mathbb{R}^{3}}(u\cdot \nabla )\theta \cdot \Delta \theta dx\leq \int_{%
\mathbb{R}^{3}}\left\vert \nabla _{h}u\right\vert \left\vert \nabla \theta
\right\vert ^{2}dx+\frac{1}{2}\int_{\mathbb{R}^{3}}\left\vert \nabla
_{h}\theta \right\vert (\left\vert \nabla u\right\vert ^{2}+\left\vert
\nabla \theta \right\vert ^{2})dx.
\end{equation*}%
Finally, we deal with the term $-\int_{\mathbb{R}^{3}}(\theta e_{3})\cdot
\Delta udx$. By integration by parts and Cauchy inequality, we have%
\begin{equation*}
-\int_{\mathbb{R}^{3}}(\theta e_{3})\cdot \Delta udx\leq \frac{1}{2}%
(\left\Vert \nabla u\right\Vert _{L^{2}}^{2}+\left\Vert \nabla \theta
\right\Vert _{L^{2}}^{2}).
\end{equation*}%
Combining the above estimates, by H\"{o}lder's inequality,
Nirenberg-Gagliardo's interpolation inequality and (\ref{eq8}), we obtain
\begin{eqnarray*}
&&\frac{1}{2}\frac{d}{dt}(\left\Vert \nabla u\right\Vert
_{L^{2}}^{2}+\left\Vert \nabla \theta \right\Vert _{L^{2}}^{2})+\left\Vert
\Delta u\right\Vert _{L^{2}}^{2}+\left\Vert \Delta \theta \right\Vert
_{L^{2}}^{2} \\
&\leq &C\int_{\mathbb{R}^{3}}(1+\left\vert \nabla _{h}u\right\vert
+\left\vert \nabla _{h}\theta \right\vert )(\left\vert \nabla u\right\vert
^{2}+\left\vert \nabla \theta \right\vert ^{2})dx \\
&\leq &C(1+\left\Vert \nabla _{h}u\right\Vert _{L^{2}}+\left\Vert \nabla
_{h}\theta \right\Vert _{L^{2}})(\left\Vert \nabla u\right\Vert
_{L^{4}}^{2}+\left\Vert \nabla \theta \right\Vert _{L^{4}}^{2}) \\
&\leq &C(1+\left\Vert \nabla _{h}u\right\Vert _{L^{2}}+\left\Vert \nabla
_{h}\theta \right\Vert _{L^{2}})(\left\Vert \nabla u\right\Vert _{L^{2}}^{%
\frac{1}{2}}\left\Vert \nabla u\right\Vert _{L^{6}}^{\frac{3}{2}}+\left\Vert
\nabla \theta \right\Vert _{L^{2}}^{\frac{1}{2}}\left\Vert \nabla \theta
\right\Vert _{L^{6}}^{\frac{3}{2}}) \\
&\leq &C(1+\left\Vert \nabla _{h}u\right\Vert _{L^{2}}+\left\Vert \nabla
_{h}\theta \right\Vert _{L^{2}})(\left\Vert \nabla u\right\Vert _{L^{2}}^{%
\frac{1}{2}}\left\Vert \nabla _{h}\nabla u\right\Vert _{L^{2}}\left\Vert
\Delta u\right\Vert _{L^{2}}^{\frac{1}{2}}+\left\Vert \nabla \theta
\right\Vert _{L^{2}}^{\frac{1}{2}}\left\Vert \nabla _{h}\nabla \theta
\right\Vert _{L^{2}}\left\Vert \Delta \theta \right\Vert _{L^{2}}^{\frac{1}{2%
}}).
\end{eqnarray*}%
Integrating this last inequality in time, we deduce that for all $\tau \in
\lbrack \Gamma ,t]$%
\begin{eqnarray}
\mathcal{J}^{2}(t) &\leq &1+\left\Vert \nabla u(\Gamma )\right\Vert
_{L^{2}}^{2}+\left\Vert \nabla \theta (\Gamma )\right\Vert _{L^{2}}^{2}+C%
\underset{\tau \in \lbrack \Gamma ,t]}{\sup }(\left\Vert \nabla _{h}u(\tau
)\right\Vert _{L^{2}}+\left\Vert \nabla _{h}\theta (\tau )\right\Vert
_{L^{2}})  \notag \\
&&\times \left( \int_{\Gamma }^{t}\left\Vert \nabla u(\tau )\right\Vert
_{L^{2}}^{2}d\tau \right) ^{\frac{1}{4}}\left( \int_{\Gamma }^{t}\left\Vert
\nabla \nabla _{h}u(\tau )\right\Vert _{L^{2}}^{2}d\tau \right) ^{\frac{1}{2}%
}\left( \int_{\Gamma }^{t}\left\Vert \Delta u(\tau )\right\Vert
_{L^{2}}^{2}d\tau \right) ^{\frac{1}{4}}  \notag \\
&&+C\underset{\tau \in \lbrack \Gamma ,t]}{\sup }(\left\Vert \nabla
_{h}u(\tau )\right\Vert _{L^{2}}+\left\Vert \nabla _{h}\theta (\tau
)\right\Vert _{L^{2}})  \notag \\
&&\times \left( \int_{\Gamma }^{t}\left\Vert \nabla \theta (\tau
)\right\Vert _{L^{2}}^{2}d\tau \right) ^{\frac{1}{4}}\left( \int_{\Gamma
}^{t}\left\Vert \nabla \nabla _{h}\theta (\tau )\right\Vert
_{L^{2}}^{2}d\tau \right) ^{\frac{1}{2}}\left( \int_{\Gamma }^{t}\left\Vert
\Delta \theta (\tau )\right\Vert _{L^{2}}^{2}d\tau \right) ^{\frac{1}{4}}
\notag \\
&\leq &1+\left\Vert \nabla u(\Gamma )\right\Vert _{L^{2}}^{2}+\left\Vert
\nabla \theta (\Gamma )\right\Vert _{L^{2}}^{2}+2C\mathcal{L}(t)\epsilon ^{%
\frac{1}{4}}\mathcal{L}(t)\mathcal{J}^{\frac{1}{2}}(t)  \notag \\
&=&1+\left\Vert \nabla u(\Gamma )\right\Vert _{L^{2}}^{2}+\left\Vert \nabla
\theta (\Gamma )\right\Vert _{L^{2}}^{2}+C\epsilon ^{\frac{1}{4}}\mathcal{L}%
^{2}(t)\mathcal{J}^{\frac{1}{2}}(t).  \label{eq26}
\end{eqnarray}%
Inserting (\ref{eq25}) into (\ref{eq26}) and taking $\epsilon $ small
enough, then it is easy to see that for all $\Gamma \leq t<T^{\ast }$, there
holds
\begin{equation*}
\mathcal{J}^{2}(t)\leq 1+\left\Vert \nabla u(\Gamma )\right\Vert
_{L^{2}}^{2}+\left\Vert \nabla \theta (\Gamma )\right\Vert
_{L^{2}}^{2}+C\epsilon ^{\frac{1}{4}}\mathcal{J}^{\frac{1}{2}}(t)+C\epsilon
^{\frac{1}{2}}\mathcal{J}^{2}(t)<\infty ,
\end{equation*}%
which proves%
\begin{equation*}
\underset{\Gamma \leq t<T^{\ast }}{\sup }(\left\Vert \nabla u(t)\right\Vert
_{L^{2}}^{2}+\left\Vert \nabla \theta (t)\right\Vert _{L^{2}}^{2})<+\infty .
\end{equation*}%
This implies us that $(u,\theta )\in L^{\infty }(0,T;H^{1}(\mathbb{R}^{3}))$%
. Thus, according to the regularity results in [2], $(u,\theta )$ is smooth
on $[0,T]$. Then we complete the proof of Theorem \ref{th1}.
\end{pf}

\section{Conclusion}

It should be noted that the condition (\ref{eq16}) is somewhat stronger than
in \cite{G1}, since it is worthy to emphasize that there are no assumptions
on the two components velocity field $(u_{1},u_{2})$. In other word, our
result demonstrates that the two components velocity field $(u_{1},u_{2})$
plays less dominant role than the one compoent velocity field does in the
regularity theory of solutions to the Boussinesq equations. In a certain
sense, our result is consistent with the numerical simulations of Alzmann et
al. in \cite{A}.

\section{Aknoledgements}

The authors are indebted to the referees for their careful reading and
valuable suggestions which improved the presentation of our paper. This work
was done while the second author was visiting the Catania University in
Italy. He would like to thank the hospitality and support of the University,
where this work was completed. This research is partially supported by
P.R.I.N. 2019. The third author wish to thank the support of "RUDN
University Program 5-100".


\begin{thebibliography}{99}
\bibitem{A} A. A. Alazman, J. P. Albert, J. L. Bona, M. Chen and J. Wu,
Comparisons between the BBM equation and a Boussinesq system,  Adv.
Differential Equations 11 (2006), 121-166.

\bibitem{BL} J. Bergh and J. L\"{o}fstr\"{o}m, Interpolation Spaces
Springer, New York, 1976.

\bibitem{BGR} S. Benbernou, S. Gala and M.A. Ragusa, On the regularity
criteria for the 3d magnetohydrodynamic equations via two components in
terms of $BMO$ space, Math. Meth. Appl. Sci. 37 (2014), 2320-2325.

\bibitem{CT0} C. Cao and E.S. Titi, Global regularity criterion for the 3D
Navier-Stokes Equations involving one entry of the velocity gradient tensor,
Arch. Rational Mech. Anal. 202 (2011), 919-932.

\bibitem{CT1} C. Cao and E.S. Titi, Regularity criteria for the three
dimensional Navier-Stokes equations, Indiana Univ. Math. J. 57 (2008),
2643-2661.

\bibitem{CN} D. Chae and H.-S. Nam, Local existence and blow-up criterion
for the Boussinesq equations, Proc. Roy. Soc. Edinburgh Sect. A 127 (1997),
935-946.

\bibitem{CKN} D. Chae, S.-K. Kim and H.-S. Nam, Local existence and blow-up
criterion of H\"{o}lder continuous solutions of the Boussinesq equations,
Nagoya Math. J. 155 (1999), 55-80.

\bibitem{CZ} J.Y. Chemin and P. Zhang, On the critical one component
regularity for the 3-D Navier-Stokes equations, Annales de l'\'{E}cole
Normale Sup\'{e}rieure 49 (2016), 133-169.

\bibitem{FO} J. Fan and T. Ozawa, Regularity criteria for the 3D
density-dependent Boussinesq equations, Nonlinearity 22 (2009), 553-568.

\bibitem{FZ} J. Fan and Y. Zhou, A note on regularity criterion for the 3D
Boussinesq system with partial viscosity, Appl. Math. Lett. 22 (2009),
802-805.

\bibitem{FQ1} D. Fang and C. Qian, Regularity criterion for 3D Navier-Stokes
equations in Besov spaces, Commun. Pure Appl. Anal. 13 (2014), 585-603.

\bibitem{FQ2} D. Fang and C. Qian, The regularity criterion for 3D
Navier-Stokes equations involving one velocity gradient component, Nonlinear
Anal. 78 (2013), 86-103.

\bibitem{GGR} S. Gala, Z. Guo and M.A. Ragusa, A remark on the regularity
criterion of Boussinesq equations with zero heat conductivity, Appl. Math.
Lett. 27 (2014), 70-73.

\bibitem{G1} S. Gala, On the regularity criterion of strong solutions to the
3D Boussinesq equations, Appl. Anal. 90 (2011), 1829-1835.

\bibitem{GR} S. Gala and M. A. Ragusa, Logarithmically improved regularity
criterion for the Boussinesq equations in Besov spaces with negative
indices, Appl. Anal. 95 (2016), 1271-1279.

\bibitem{He} C. He, Regularity for solutions to the Navier-Stokes equations
with one velocity component regular, Electr. J. Differential Equations, 2002
(2002), 1-13.

\bibitem{H} R. A. Hunt, On L(p,q) spaces, Enseign. Math., 12 (1966), 249-276.

\bibitem{IM} N. Ishimura and H. Morimoto, Remarks on the blow-up criterion
for 3D Boussinesq equations, Math. Models Methods Appl. Sci. 9 (1999),
1323-1332.

\bibitem{JZ1} X.J. Jia and Y. Zhou, Remarks on regularity criteria for the
Navier-Stokes equations via one velocity component, Nonlinear Anal. Real
World Appl. 15 (2014), 239-245.

\bibitem{KZ} I. Kukavica and M. Ziane, One component regularity for the
Navier-Stokes equations, Nonlinearity 19 (2006) 453-469.

\bibitem{NP} J. Neustupa and P. Penel, Regularity of a suitable weak
solution to the Navier-Stokes equations as a consequence of regularity of
one velocity component, Applied Nonlinear Analysis, edited by Sequeira et
al. Kluwer Academic / Plenum Publishers, New York, 1999.

\bibitem{O} R. O'Neil, Convolution operators and $L^{p,q}$ spaces, Duke
Math. J. 30 (1963), 129-142.

\bibitem{PP} P. Penel and M. Pokorn\'{y}, On anisotropic regularity criteria
for the solutions to 3D Navier-Stokes equations, J. Math. Fluid Mech. 13
(2011), 341-353.

\bibitem{Q} C. Qian , A generalized regularity criterion for 3d
Navier-Stokes equations in terms of one velocity component, J. Differ. Equ.
260 (2016), 3477-3494 .

\bibitem{QYW} Y. Qin, X. Yang, Y. Wang, X. Liu, Blow-up criteria of smooth
solutions to the 3D Boussinesq equations, Math. Methods Appl. Sci. 35
(2012), 278-285.

\bibitem{QDY} H. Qiu, Y. Du and Z. Yao, A blow-up criterion for 3D
Boussinesq equations in Besov spaces, Nonlinear Anal. 73 (2010), 806-815.

\bibitem{Tri} H. Triebel, Theory of Function Spaces, Birkh\"{a}user Verlag,
Basel, Boston, 1983.

\bibitem{X} Z. Xiang, The regularity criterion of the weak solution to the
3D viscous Boussinesq equations in Besov spaces, Math. Methods Appl. Sci. 34
(2011), 360-372.

\bibitem{XZZ} F. Xu, Q. Zhang and X. Zheng, Regularity criteria of the 3D
Boussinesq equations in the Morrey-Campanato space, Acta Appl. Math. 121
(2012), 231-240.

\bibitem{Ye3} Z. Ye, Blow-up criterion of smooth solutions for the
Boussinesq equations, Nonlinear Anal. 110 (2014), 97-103.

\bibitem{Ye} Z. Ye, A logarithmically improved regularity criterion of
smooth solutions for the 3D Boussinesq equations, Osaka J. Math. 53 (2016),
417-423.

\bibitem{Ye1} Z.\ Ye, Remarks on the regularity criterion to the 3D
Navier-Stokes equations via one velocity component, Differential Integral
Equation, 29 (2016), 957-976.

\bibitem{Ye2} Z.\ Ye, Remarks on the regularity criterion to the
Navier-Stokes equations via the gradient of one velocity component, J. Mah.
Anal. Appl. 435 (2016), 1623-1633.

\bibitem{Zhe} X. Zheng, A regularity criterion for the tridimensional
Navier-Stokes equations in terms of one veloicty component, J. Differential
Equations, 256 (2014), 283-309.

\bibitem{Z1} Y. Zhou, A new regularity criterion for weak solutions to the
Navier-Stokes equations, J. Math. Pures Appl. 84 (2005), 1496-1514.

\bibitem{ZP1} Y. Zhou and M. Pokorn\'{y}, On the regularity of the solutions
of the Navier-Stokes equations via one velocity component, Nonlinearity 23
(2010), 1097-1107.
\end{thebibliography}
\end{document}